\documentclass[11pt]{article}

\input xy
\xyoption{knot}

\usepackage{fullpage, amsfonts, amsmath, amsthm,stmaryrd}
\usepackage{fullpage, amsfonts, amsmath, amsthm,stmaryrd}
\usepackage[all]{xy}
\usepackage{epsf}
\usepackage{graphicx}
\usepackage{amssymb}
\usepackage{amscd}
\usepackage{color}

\newtheorem{thm}{Theorem}[section]

\newtheorem{prop}[thm]{Proposition}
\newtheorem{defin}[thm]{Definition}

\def\c{{\mathcal C}}
\def\f{{\mathcal F}}

\def\D{{\mathcal D}}
\def\F{{\mathbb F}}
\def\H{{\mathbb H}}
\def\L{{\mathbb L}}

\def\P{{\mathcal P}}

\def\R{{\mathbb R}}
\def\T{{\mathcal T}}
\def\Z{{\mathbb Z}}

\def\del{{\partial}}

\def\mod{{\textup{mod} \;}}

\begin{document}

\title{A surgery triangle for lattice cohomology}

\author{Joshua Greene}

\date{\today}

\maketitle

\begin{abstract}
Lattice cohomology, defined by N\'emethi in \cite{Nlattice}, is an invariant of negative definite plumbed 3-manifolds which conjecturally computes the Heegaard Floer homology $HF^+$.  We prove a surgery exact triangle for the lattice cohomology analogous to the one for $HF^+$.  This is a step towards comparing these two invariants.
\end{abstract}

\section{Introduction.}

In this note, we study a model for the Heegaard Floer homology $HF^+$ of a {\em negative definite plumbed manifold}.  Heegaard Floer homology is a powerful invariant in low-dimensional topology, and the hat version $\widehat{HF}$ is now known to be algorithmically computable for an arbitrary $3$-manifold \cite{SW}.  However, the same is unknown for the structurally richer invariant $HF^+$, and negative definite plumbed manifolds are one family for which calculation seems feasible.  Moreover, calculations in this case are useful in applications to knot concordance \cite{GJ}, unknotting number \cite{OSunknotting}, and smoothings of rational surface singularities \cite{Stipsicz}.  In this setting, N\'emethi has given an elegant construction of an algebraic invariant called the {\em lattice cohomology} $\H^+$ which conjecturally computes $HF^+$ \cite{Nlattice}, and this is the invariant we study here.

To begin, we describe the family of manifolds at hand.  Let $G = (V,E)$ denote a finite planar graph with an integer weight $m(v)$ for each vertex $v \in V$ and a sign $\pm$ for each edge $e \in E$.  Associated to $G$ is a compact 4-manifold with boundary $X(G)$.  To describe it, center a round unknot at each vertex of $G$ with framing $m(v)$, and introduce a right-hand (resp. left-hand) clasp between neighboring unknots for each positive (resp. negative) edge $e \in E$.  Let $\L$ denote the resulting framed link in $S^3 = \del D^4$, and let $X(G)$ denote the result of attaching 2-handles to $D^4$ along $\L$.  When $G$ is {\em acyclic} (i.e. a disjoint union of trees), $W(G)$ admits an alternative description: form the disk bundle over $S^2$ with Euler number $m(v)$ for each $v \in V$, and plumb together the bundles of neighboring vertices for each edge $e \in E$.  Note that the signs on the edges of $G$ become immaterial in this case.  We let $Y(G)$ denote the oriented boundary of $X(G)$.

The graph $G$ defines a free $\Z$-module $L = L(G)$ which is freely generated by classes $E_v$, $v \in V$.  There is a symmetric bilinear form on this module which is defined by setting $(E_v,E_w)$ equal to the signed number of edges between $v$ and $w$, when $v \ne w$, and $(E_v,E_v)=m(v)$.  The quadratic module $(L(G),(\cdot,\cdot))$ defined in this way is isomorphic to the module $H_2(X(G);\Z)$ equipped with its intersection pairing.  The isomorphism is set up by identifying a class $E_v$ with $[S_v]$, where $S_v$ is the sphere obtained by taking the union of a Seifert disk for $K_v$ with the core of the corresponding handle attachment.  In the case of the plumbing description, $S_v$ may be alternately viewed as the base sphere in the disk bundle associated to $v \in V$.

\begin{defin}

A graph $G$ is {\em negative-definite} if it is acyclic and the pairing $(\cdot,\cdot)$ is negative definite, and the resulting space $Y(G)$ is a {\em negative definite plumbed manifold}.

\end{defin}

\noindent The spaces $Y(G)$ associated to negative definite graphs are precisely the boundaries of {\em normal surface singularities} \cite{Grauert, Mumford}; for a nice summary see \cite[esp. pp. 282-283]{Nlectures}.

Ozsv\'ath-Szab\'o undertook the determination of $HF^+$ of a negative definite plumbed manifold under an additional assumption on the graph $G$ \cite{OSplumbed}.  Namely, they call a vertex $v \in V$ {\em bad} if $m(v) + d(v) > 0$; here $d(v)$ denotes the number of edges incident $v$.  Their result is a calculable description of the invariant $HF^+(-Y(G))$ under the assumption that $G$ has at most one bad vertex.  Additionally, they describe the ``even half" of this invariant $HF^+_{even}(-Y(G))$ when $G$ has two bad vertices.  Using their result, Ozsv\'ath-Szab\'o calculated the invariant for many spaces of interest, and gave a useful formula for the {\em correction terms} of such a space.

However, it remained a problem to find a suitable generalization of their algorithm which would apply to {\em any} negative definite plumbed 3-manifold $Y(G)$.  This was taken up by N\'emethi, who proposed the lattice cohomology $\H^+(G)$.  He proved that the lattice cohomology is a diffeomorphism invariant of $Y(G)$ (independent of the graph $G$), reduces to Ozsv\'ath-Szab\'o's model in the known case, and also agrees with $HF^+$ for the boundaries of rational and elliptic surface singularities \cite{Nplumbed, Nlattice}.  As a result, N\'emethi was able to show that rational surface singularities are L-spaces; and his conjecture $\H^+(G) \cong HF^+(-Y(G))$ would imply that these are {\em all} the L-spaces amongst negative definite plumbed manifolds.  At present there is no purely topological characterization of L-spaces.

An essential feature of the Floer homology groups is that they obey a surgery exact triangle: if $K \subset Y$ is a framed knot, then there is a long exact sequence \[ \cdots \to HF^+(Y) \to HF^+(Y_0(K)) \to HF^+(Y_1(K)) \to HF^+(Y) \to \cdots; \] here $Y_r(K)$ denotes the result of $r$-framed surgery on $K$.  As a first step towards proving N\'emethi's conjecture, it is desirable to know that the lattice cohomology obeys an analogous exact triangle.  Then the hope is to define a map $HF^+ \to \H^+$ and establish that this map is an isomorphism by an induction using the five-lemma.

The purpose of this note is to prove Theorem \ref{thm:triangle}, which establishes the existence of the surgery exact triangle for the lattice cohomology over the coefficient ring $\Z / 2 \Z$.  Moreover, we describe a version of the lattice cohomology which applies to the family of manifolds $Y(G)$ without any additional assumption on the graph $G$, and show that the triangle holds in this setting.  Indeed, the proof of the triangle is a matter of fairly straightforward (if somewhat involved) algebra once we have suitably defined the invariant $\H^+$ for arbitrary $Y(G)$.  We fully expect an analogue of Theorem \ref{thm:triangle} to hold with $\Z$ coefficients as well, a matter we intend to take up in future work.  What is less clear at this time is how to define the desired map $HF^+ \to \H^+$, even in the primary case of a negative definite plumbed manifold.  This is a subject for further study.

Finally, and what is perhaps most intriguing, the lattice cohomology comes equipped with an additional integer grading which has no obvious counterpart in the Floer homology.  If $L \subset S^3$ is a link, then the double cover of $S^3$ branched along $L$ takes the form $Y(G)$ for some $G$, and conversely all negative-definite plumbed manifolds arise in this way.  For branched double-covers there is a conjectural $\delta$-grading described in \cite{Baldwin} and \cite[Conjecture 8.1]{Greene}.  Indeed, these two quantities agree (up to an overall shift) in the limited domain where they have both been computed.  In light of the way we expect the $\delta$-grading to arise from \cite{Baldwin}, it is compelling to search for a more direct relationship between the Khovanov homology and lattice cohomology.

\subsection*{Acknowledgments.}

Thanks to my advisor Zolt\'an Szab\'o for encouraging me to pursue this project, and to Andr\'as N\'emethi for helpful correspondence.

\section{The lattice cohomology.}\label{ss:lattice cohomology}

The pairing on $L(G)$ gives rise to a quadratic form $q(x) = -\frac{1}{8} (x,x)$, and extends to a pairing on $L \otimes \R$.  A vector $K \in \text{Hom}(L,\Z)$ is {\em characteristic} if $(K,x) \equiv (x,x) \; (\mod 2)$ for all $x \in L(G)$; and the set of all characteristic vectors is denoted $Char(G)$.  Define \[ \c^+(G) = \text{Hom}_\F(Char(G) \times \P(V),\T^+_0).\]  Here $\F$ denotes either $\Z$ or $\Z / 2 \Z$; $\P(V)$, the power set of $V$; and $\T^+_0$, the $\F[U]$-module $\F[U,U^{-1}] / U \cdot \F[U]$.  The latter is a graded $\F$-module, with grading given by $\text{gr}(U^{-d}) = 2d$, and makes $\c^+(G)$ into an $\F[U]$-module.  The group $\c^+(G)$ has an additional grading, where $\phi \in \c^+(G)$ is homogeneous of degree $s$ if it is supported on pairs $(K,S)$ with $|S|=s$.

Consider the affine lattice $K + 2L \subset L \otimes \R$.  Associated to a pair $(K',S) \in (K + 2L) \times \P(V)$ is an $s$-dimensional cube in $L \otimes \R$ with vertex set $\{ K' + 2 \sum_{j \in T} E_j \; | \; T \subset S\}$. In this way we obtain a cubical decomposition of $L \otimes \R$, and in turn a chain complex $(C_s,\del)$.  Under the identification between pairs $(K,S)$ and cubes, we can express the boundary operator as \[ \del(K,S) = \sum_{(K',S')} \epsilon(K',S') \cdot (K',S'),\] where $\epsilon(K',S') = \pm 1$ for pairs of the form $(K,S-w)$ or $(K+2E_w,S-w)$, $w \in S$, and $\epsilon(K',S') = 0$ otherwise.  We extend the quadratic form $q$ to a function on $Char(G) \times \P(V)$ by setting \[ q(K,S) = \max \{ q(K + 2 \sum_{j \in T} E_j) \; | \; T \subset S\}. \]  Lastly, we define a differential $\delta$ on the group $\c^+(G)$ by setting \[ \delta(\phi)(K,S) = \sum_{(K',S')} \epsilon(K',S') \cdot U^{q(K,S)-q(K',S')} \cdot \phi(K',S').\]  It is straightforward to verify that $\delta^2=0$, making use of the fact that $\del^2=0$ \cite[Lemma 3.1.5]{Nlattice}.  The following definition is in essence \cite[Definition 3.2.5]{Nlattice}.

\begin{defin}\label{def: H^+}
The homology of $\c^+(G)$, regarded as a graded $\F[U]$-module, is the {\em lattice cohomology} of $G$.  It is denoted by $\H^+(G)$.
\end{defin}

\noindent Observe that this definition makes sense for {\em any} finite graph $G$, with no assumption on planarity; however, it is unclear what significance it has in this more general setting.  In the case that $G$ is negative definite, there is a more topological description of the invariant $\H^+(G)$ \cite[Definitions 3.1.11]{Nlattice}, though we do not use it here.  This is the invariant we wish to compare to $HF^+(-Y(G))$.  In fact, \cite[Conjecture 5.2.4]{Nlattice} proposes in a precise form the conjecture that $\H^+(G) \cong HF^+(-Y(G))$ in the case that $G$ is negative definite.

\section{Statement of the triangle.}\label{ss:triangle}

The invariant $HF^+$ obeys an exact triangle.  For the case of interest, this triangle takes the following form.  Select $v \in V$, and let $K$ denote a meridian for the link component $K_v$.  Then $0$-surgery on $K$ is $Y(G-v)$ and $(-1)$-surgery on $K$ is $Y(G_{+1}(v))$, where $G-v$ denotes the graph gotten by deleting $v$ and all its incident edges, and $G_{+1}(v)$ denotes the graph gotten by increasing $m(v)$ by one.  Thus the triangle reads: \[\cdots \to HF^+(-Y(G_{+1}(v))) \to HF^+(-Y(G)) \to HF^+(-Y(G-v)) \to HF^+(-Y(G_{+1}(v))) \to \cdots \]  Our purpose here is to prove an analogous result for the lattice cohomology.

\begin{thm}\label{thm:triangle}

Let $\F$ denote the field $\Z / 2 \Z$.  There is a short exact sequence of complexes \[ 0 \to \c^+(G_{+1}(v)) \stackrel{A}\to \c^+(G) \stackrel{B}\to \c^+(G-v) \to 0\] which gives rise to an $\F[U]$-equivariant exact triangle of lattice cohomology groups \[ \cdots \to \H^+(G_{+1}(v)) \to \H^+(G) \to \H^+(G-v) \to \H^+(G_{+1}(v)) \to \cdots \]

\end{thm}

\noindent The rest of this note is devoted to the proof of Theorem \ref{thm:triangle}, which we break into several pieces.

\section{Definition of the maps appearing in the short exact sequence.}\label{ss:A and B}

There is a canonical identification between $L(G-v)$ and the sublattices of $L(G)$ and $L(G_{+1}(v))$ spanned by the classes $E_w, w \ne v$.  We denote a characteristic vector for $G$ or $G_{+1}(v)$ by a pair $(K,t)$, where $K$ denotes the restriction to $L(G-v)$ and $t$ is the pairing of the vector with $E_v$.  Denote by $(\cdot,\cdot)'$ the pairing on $L(G_{+1}(v))$ and by $q'$ the associated function on $Char(G_{+1}(v)) \times \P(V)$.

The definition of the maps $A$ and $B$ on the $0^{th}$ level of the short exact sequence is taken from \cite{OSplumbed}.  We set \begin{equation}\label{eqn:A0}
A(\phi)(K,t) = \sum_{i = -\infty}^{+ \infty} U^{i(i+1)/2} \cdot \phi(K,t+2i+1)
\end{equation} and

\begin{equation}\label{eqn:B0}
B(\phi)(K) = \sum_{i=-\infty}^{+ \infty} \phi(K,m(v)+2i).
\end{equation} (We suppress $S = \emptyset$ from the notation.) Observe that $A = \mathbb{A}^+ \circ R$ and $B = \mathbb{B}^+$, where the maps $\mathbb{A}^+$ and $\mathbb{B}^+$ appear just before Lemma 2.9 and $R$ appears in Proposition 2.5 of \cite{OSplumbed}.

Let us attempt to extend this definition to a pair of $\F[U]$-equivariant chain maps $A: \c^+(G_{+1}(v)) \to \c^+(G)$ and $B: \c^+(G) \to \c^+(G-v)$.  We focus on the definition of $A$ first.  The condition $\delta A = A \delta$ and a staightforward induction on $|S|$ shows that the value of $A(\phi)((K,t),S)$ is an $\F[U]$-linear combination of terms $\phi((K,t+2i+1),S)$.  Moreover, the coefficient on $\phi((K,t+2i+1),S)$ is a monomial $U^{c(i,(K,t),S)}$, for some $c(i,(K,t),S) \geq 0$.  We determine the value of this exponent by comparing the coefficients on $\phi((K,t+2i+1),S-w)$ on both sides of the identity $\delta A(\phi)((K,t),S) = A \delta(\phi)((K,t),S)$, for $w \in S$.  An induction on $|S|$ shows that this value is uniquely determined by the expression

\begin{equation}\label{eqn:c(i)}
\begin{split}
c(i,(K,t),S) =  [q((K,t),S) - & q((K,t))] - [q'((K,t+2i+1),S) - q'(K,t+2i+1)] \\ & + i(i+1)/2.
\end{split}
\end{equation}  However, both sides of the identity $\delta A(\phi)((K,t),S) = A \delta(\phi)((K,t),S)$ also involve terms of the form $\phi((K,t+2i+1)+2E_w,S-w)$.  It stands to check that the stated definition of $c(i,(K,t),S)$ makes the coefficients on these terms agree as well.  This is easily confirmed for the case when $w \ne v$.  In the case $w=v$, there is one small wrinkle.  Namely, we must be careful to recognize that the class $(K,t+2i-1)+2E_v \in Char(G_{+1}(v))$ is the one whose evaluation on $E_w$ with respect to the pairing $(\cdot,\cdot)'$ is given by $(K+2E_v,E_w)$ when $w \ne v$ and $(K+2E_v,E_v)+2i+1$ when $w = v$.  This is due to the fact that $(E_v,E_v)' = (E_v,E_v)+1$.  Having observed this, the case $w=v$ follows as well.

Next, we need to verify that $c(i,(K,t),S) \geq 0$, so that $A$ is a bona fide $\F[U]$-module map, rather than just an $\F[U,U^{-1}]$-module map.  The value of $q((K,t),S) - q(K,t)$ is the maximum of the values $q((K,t) +2E_T) - q(K,t)$, where $E_T := \sum_{j \in T} E_j, T \subset S.$
It is easy to check that

\begin{eqnarray*}
q'((K,t+2i+1) +2E_T) - q'(K,t+2i+1) = \left\{
\begin{array}{cl}
q((K,t) +2E_T) - q(K,t), & \quad \mbox{if $v \notin T$}; \cr
& \cr
q((K,t) +2E_T) - q(K,t)-i-1, & \quad \mbox{if $v \in T$}. \cr
\end{array}
\right.
\end{eqnarray*}  It follows that if $v \notin S$, then \[ q'((K,t+2i+1),S) - q'(K,t+2i+1) = q((K,t),S) - q(K,t),\] and if $v \in S$, then \[ q'((K,t+2i+1),S) - q'(K,t+2i+1) = \max \{q((K,t),S-v),q((K,t)+2E_v,S-v)-i-1 \} - q(K,t).\]  Denote by $r((K,t),S)$ the difference $q((K,t),S-v)-q((K,t)+2E_v,S-v)$ when $v \in S$.  We conclude that if $v \notin S$, then

\begin{equation}\label{eqn:c(i)1}
c(i,(K,t),S) = i(i+1)/2;
\end{equation} and if $v \in S$, then

\begin{eqnarray}\label{eqn:c(i)2}
c(i,(K,t),S) = \left\{
\begin{array}{cl}
i(i+1)/2, & \quad \mbox{if $r((K,t),S) \geq \max \{0, -(i+1) \}$}; \cr
& \cr
(i+1)(i+2)/2, & \quad \mbox{if $r((K,t),S) \leq \min \{ 0, -(i+1) \}$}; \cr
& \cr
(i+1)(i+2)/2 + r((K,t),S), & \quad \mbox{if $0 \leq r((K,t),S) \leq -(i+1)$}; \cr
& \cr
i(i+1)/2 - r((K,t),S), & \quad \mbox{if $-(i+1) \leq r((K,t),S) \leq 0$}. \cr
\end{array}
\right.
\end{eqnarray}

\noindent Observe that in any case, $c(i,(K,t),S) \geq 0$, so that $A$ is indeed an $\F[U]$-module map.  Moreover, we can tell exactly when this value vanishes.  Fix $(K,t)$ and $S$, and abbreviate $r = r((K,t),S)$.  If $v \notin S$, then $c(i,(K,t),S) = 0$ iff $i = 0$ or $-1$.  If $v \in S$, then $c(i,(K,t),S) = 0$ iff $i=0$ and $r \geq 0$; $i=-1$; or $i=-2$ and $r \leq 0$.

The extension of $B$ to a chain map $\c^+(G) \to \c^+(G-v)$ proceeds similarly, but is simpler.  The result is summarized as follows.

\begin{prop}\label{prop:AB}

There is a unique way to extend the maps $A$ and $B$ in Equations (\ref{eqn:A0}) and (\ref{eqn:B0}) to $\F[U]$-equivariant chain maps $A: \c^+(G_{+1}(v)) \to \c^+(G)$ and $B: \c^+(G) \to \c^+(G-v)$.  They are defined by setting

\begin{equation}\label{eqn:A}
A(\phi)((K,t),S) = \sum_{i=-\infty}^{+ \infty} U^{c(i,(K,t),S)} \cdot \phi((K,t+2i+1),S),
\end{equation} where $c(i,(K,t),S)$ is defined by Equation (\ref{eqn:c(i)}), and

\begin{equation}\label{eqn:B}
B(\phi)((K,t),S) = \sum_{i=-\infty}^{+ \infty} \phi((K,m(v)+2i+1),S).
\end{equation}

\end{prop}

\section{Proof of exactness.}

\subsection{$A$ is injective.}\label{ss:A inj}

Select a non-zero element $\phi \in \c^+(G)$.  There is a filtration $\f$ on $\T^+_0$ induced by the grading $\text{gr}(U^{-d})=2d$, and where by convention $\f(0) = -\infty$.  Choose a pair $((K,t),S)$ such that (a) $\f(\phi((K,t),S))$ is maximal and (b) $t$ is minimal subject to condition (a).  Set $r = r((K,t-1),S)$.  Let us suppose first that $r \geq 0$.  Consider \[ \f(A(\phi)((K,t-1),S)) = \f \left( \sum_{i=-\infty}^{+ \infty} U^{c(i,(K,t-1),S)} \cdot \phi((K,t+2i),S) \right). \]  The remarks following Equation (\ref{eqn:c(i)2}) imply that the $i=0$ term is just $\phi((K,t),S)$.  Those terms with $i > 0$ have lower $\f$-grading by condition (a) and the fact that $c(i,(K,t-1),S) > 0$; and those with $i < 0$ have lower $\f$-grading by condition (b).  Consequently, the value of $\f(A(\phi)((K,t-1),S))$ simply reduces to $\f(\phi((K,t),S)$).  It follows that $A(\phi)((K,t-1),S)) \ne 0$, since we assumed $\phi \ne 0$.  If instead $r \leq -1$, we repeat the same argument, only we choose $t$ minimal so that $\f(\phi(K,t-2),S)$ is maximal.  In this case it is the $i=-1$ term with maximal $\f$-grading.  The conclusion is that $A(\phi) \ne 0$ in any case, so $A$ is injective.



\subsection{$B$ is surjective.}\label{ss:B surj}

Denote by $(K,S)^\vee$ the $\text{Hom}$-dual to the pair $(K,S)$.  Thus $\c^+(G-v)$ is freely generated as an $\F$-module by the maps $U^{-m} \cdot (K,S)^\vee$, $m \geq 0$.  It is clear that \[ B(U^{-m} \cdot ((K,t),S)^\vee) = U^{-m} \cdot (K,S)^\vee \] for any $(K,S)$ and $t \equiv m(v) \; (\mod 2)$.  Hence $B$ is surjective.

\subsection{$\text{im} \; A \subset \text{ker} \; B$.}\label{ss:im A and ker B}
We compute

\begin{eqnarray*}
B \circ A (\phi) (K,S) &=& \sum_{i = - \infty}^{+ \infty} \sum_{j = - \infty}^{+ \infty} U^{c(j,(K,m(v)+2i+1),S)} \cdot \phi((K,m(v)+2i+2j+1),S) \\
&=& \sum_{k = - \infty}^{+ \infty} \sum_{j = - \infty}^{+ \infty} U^{c(j,(K,m(v)+2(k-j)+1),S)} \cdot \phi((K,m(v)+2k+1),S). \\
\end{eqnarray*}  Since $S \subset V(G-v)$, we have by Equation (\ref{eqn:c(i)1}) that $c(j,(K,m(v)+2(k-j)+1),S) = j(j+1)/2$ for all $k$ and $j$.  For fixed $k$, the multiplier $U^{j(j+1)/2}$ annihilates $\phi((K,m(v)+2k+1),S)$ for $|j| \gg 0$.  Those which do not annihilate this term cancel in pairs since $j(j+1)/2$ is symmetric in $j$ and $-1-j$. Hence $B \circ A (\phi) = 0$.

\subsection{$\text{ker} \; B = \D$.}\label{ss:ker B and D}

We define $\D_1 \subset \c^+(G)$ to be the submodule freely generated over $\F$ by the mappings  \[U^{-m} \cdot ((K,t),S)^\vee, \quad v \in S, \] $\D_2 \subset \c^+(G)$ to be the submodule freely generated over $\F$ by \[ U^{-m} \cdot [((K,t),S)^\vee + ((K,t+2),S)^\vee], \quad v \notin S, \]  and $\D = \D_1 \oplus \D_2$.  It is easy to see that $\phi \in \D$ iff in the expansion of $\phi$ with respect to the basis $\{ U^{-m} \cdot ((K,t),S)^\vee \}$ of $\c^+(G)$, there are an even number of terms $U^{-m} \cdot ((K,t+2i),S)^\vee$ as $i$ varies and $m,(K,t)$, and $v \notin S$ remain fixed.  On the other hand, $\text{ker} \; B$ clearly fulfills this description as well.  Hence $\text{ker} \; B = \D$.

\subsection{$\D \subset \text{im} \; A$.}\label{ss:D and im A}

Filter $\c^+(G)$ by letting $\c^{(m)}$ consist of those maps whose image lies in the $m^{th}$ filtered piece of $\T_0^+$, and set $\D_1^{(m)} = \D_1 \cap \c^{(m)}, \D_2^{(m)} = \D_2 \cap \c^{(m)}$, and $\D^{(m)} = \D \cap \c^{(m)}$.  We prove by induction on $m$ that $\D^{(m)} \subset \text{im} \; A$.

Fix a pair $((K,t),S)$ with $v \notin S$, and consider the map $((K,t+1),S)^\vee \in \c^+(G_{+1}(v))$.  Observe that in this case Equation (\ref{eqn:c(i)1}) computes the value of $c(i,(K,t),S)$.  From this it is clear that $A(((K,t+1),S)^\vee)$ vanishes on all pairs $((K',t'),S')$ not of the form $((K,t+2i),S)$, and for this pair it takes the value $1$ if $i=0$ or $1$, and $0$ otherwise.  In other words, $A(((K,t+1),S)^\vee) = ((K,t),S)^\vee + ((K,t+2),S)^\vee.$  Hence $\D_2^{(0)} \subset \text{im} \; A$.

Next, fix a pair $((K,t),S)$ with $v \in S$, and consider the map $((K,t+2i+1),S)^\vee \in \c^+(G_{+1}(v))$.  Once again, the image of this map under $A$ vanishes on all pairs not of the form $((K,t+2j),S)$, and for a pair of this form it takes the value $U^{c(i-j,(K,t+2j),S)}$.  This value vanishes unless $c(i-j,(K,t+2j),S) = 0$, and we can determine when this occurs by way of the remarks following Equation (\ref{eqn:c(i)2}).  First, observe that $r((K,t+2j),S) = r((K,t),S) +j$ (recall the definition just before Equation (\ref{eqn:c(i)1})).  Therefore, in working with the family of pairs $((K,t+2j),S)$ as $j$ varies, we may assume that $K$ is chosen so that $r((K,t),S) = 0$.  With this choice made, the remarks following Equation (\ref{eqn:c(i)2}) implies that $c(i-j,(K,t+2j),S) = 0$ iff (a) $j = i \geq 0$; (b) $j = i+1$; or (c) $j = i+2, i \leq -2$.  It follows that

\begin{eqnarray*}
A(((K,t+2i+1),S)^\vee) = \left\{
\begin{array}{cl}
((K,t+2i),S)^\vee + ((K,t+2i+2),S)^\vee, & \quad i \geq 0; \cr
& \cr
((K,t),S)^\vee, & \quad i = -1; \cr
& \cr
((K,t+2i+2),S)^\vee + ((K,t+2i+4),S)^\vee, & \quad i \leq -2. \cr
\end{array}
\right.
\end{eqnarray*}  By taking linear combinations of these terms, we see at once that $((K,t+2i),S)^\vee \in \text{im} \; A$ for all $i$. Hence $\D_1^{(0)} \subset \text{im} \; A$, and so $\D^{(0)} \subset \text{im} \; A$.

Now suppose that we have shown that $\D^{(m-1)} \subset \text{im} \; A$ for some $m \geq 1$, and choose $\phi \in \D^{(m)} - \D^{(m-1)}$.  Then $U \cdot \phi \in \D^{(m-1)} \subset \text{im} \; A$: say $A(\psi) = U \cdot \phi$.  Then $A( U^{-1} \cdot \psi) = \phi + \rho$ for some $\rho \in \c^{(0)}$.  Since $\text{im} \; A \subset \D$ by $\S$\ref{ss:im A and ker B} and $\S$\ref{ss:ker B and D}, we have $A( U^{-1} \cdot \psi) \in \D$, and so $\rho = A( U^{-1} \cdot \psi) - \phi \in \D^{(0)}$.  On the other hand, we have just seen that $\D^{(0)} \subset \text{im} \; A$.  Hence $\phi = A( U^{-1} \cdot \psi) - \rho \in \text{im} \; A$, completing the induction step.  In conclusion, $\D \subset \text{im} \; A$.

\vspace{1cm}

Collecting the results of the previous subsections, we conclude that the chain maps $A$ and $B$ appearing in Proposition \ref{prop:AB} provide the desired maps for the short exact sequence in Theorem \ref{thm:triangle}.  The resulting exact triangle of lattice cohomology groups follows directly.

\end{document}